\documentclass{article}

\usepackage{amsmath,amsthm,amssymb}

\newtheorem{prop}{Proposition}[section]
\newtheorem{cor}{Corollary}[section]

\newcommand{\eps}{\varepsilon}

\newcommand{\ph}{\varphi}

\newcommand{\const}{\mathrm{const}}

\newcommand{\R}{{\mathbb{R}}}
\newcommand{\T}{{\mathbb{T}}}

\newcommand{\calA}{\mathcal{A}}

\newcommand{\calF}{\mathcal{F}}

\newcommand{\calS}{\mathcal{S}}
\newcommand{\calT}{\mathcal{T}}

\begin{document}

\title{Another Billiard Problem}
\author{Sergey Bolotin and Dmitry Treschev\\
Steklov Mathematical Institute of Russian Academy of Sciences}

\date{ }
\maketitle

\begin{abstract}
Let $(M,g)$ be a Riemannian manifold, $\Omega\subset M$ a domain with boundary $\Gamma$,
and $\phi$ a smooth function such that $\phi|_\Omega > 0$, $\ph|_\Gamma = 0$, and $\nabla\phi|_\Gamma\ne 0$.
We study the geodesic flow of the metric $G=g/\phi$. The $G$-distance from any point of $\Omega$ to $\Gamma$ is finite,
hence the geodesic flow is incomplete. Regularization of the flow in a neighborhood of $\Gamma$ establishes a natural reflection law from $\Gamma$. This leads to a certain (not quite standard) billiard problem in $\Omega$.

\end{abstract}


\section{Introduction}

Billiard systems were introduced by Birkhoff \cite{Birk}. They became an important class of models in Hamiltonian dynamics.
Later several generalizations appeared: outer billiards \cite{Moser}, magnetic billiards \cite{magn}, projective billiards \cite{proj} and many others, see \cite{Tab}.
Another addition to this list are billiard-like systems motivated by the problem of water waves in a shallow reservoir.
These systems were studied in the papers \cite{Dob1,Dob2,Naz1,Naz2,Naz3}.

Let $\Omega$ be an open domain in $\R^n$ with smooth boundary $\Gamma=\partial\Omega$. The Birkhoff billiard in $\Omega$ is defined as follows.
A particle moves in $\Omega$ with constant unit velocity $v$, and when colliding with the boundary at a point $q\in \Gamma$,
it is reflected so that the energy  and the tangent momentum\footnote{$n(q)$ is the unit normal vector to $\Gamma$.} $p=v-\langle v,n(q)\rangle n(q)$ are conserved.
Trajectories of the billiard are determined by the symplectic billiard map $(q,p)\mapsto (q_+,p_+)$,
where $q_+$ is the next collision point, and $p_+$ the corresponding tangent momentum.
The billiard map is everywhere defined and smooth on $\{(q,p)\in T^*\Gamma:|p|<1\}$ when the boundary is strictly convex.

What makes Birkhoff billiard important is that it appears in the short wave approximation to the wave equation in $\Omega$:
$$
u_{tt}=\Delta u,\qquad u|_\Gamma=0.
$$
For example, a  Lagrangian invariant torus of the billiard map gives a series of quasi-eigenfunctions of the Laplace operator.

The problem studied in the present paper is motivated by the  short wave approximation for the wave equation in $\Omega$ degenerating on the boundary:
\begin{equation}
\label{eq:wave}
  u_{tt}=\nabla_g\cdot (\phi\nabla_g u),
\end{equation}
Here $\phi$ is a smooth function such that $\phi>0$ in $\Omega$, $\phi|_\Gamma=0$, and $d\phi|_\Gamma\ne 0$.
The divergence   and the gradient $\nabla_g$  are determined by some Riemannian metric $g$, not necessarily Euclidean.
To obtain solutions of (\ref{eq:wave}) one has to know eigenvalues and eigenfunctions of the differential
operator\footnote{The operator degenerates on the boundary, and instead of the boundary conditions its Friedrichs self-adjoint extension should be used.}
in (\ref{eq:wave}).

This problem was studied in \cite{Dob1,Dob2,Naz1,Naz2,Naz3}.
To apply Maslov's canonical operator method in the short wave approximation, one needs to find Lagrangian submanifolds
in the level set of the Hamiltonian\footnote{We use the same notation for the norms defined by the metric $g$ in $T_x\Omega$
and $T_x^*\Omega$.}  $H(x,p_x)=\phi(x)\|p_x\|_g^2$ on $T^*\Omega$.
Such submanifolds are invariant under the flow of the Hamiltonian system
$$
   \dot x = \partial_{p_x} H,\quad
  \dot p_x = -\partial_x H.
$$
The flow is incomplete: solutions hit the boundary in finite time. In \cite{Naz2} it was shown that the flow is regularizable.
In this paper we introduce a simpler regularization. The regularized configuration space is the double of $\Omega$,
and the regularized geodesic flow is  complete.

The dynamical system obtained in this way is analogous to the traditional billiard systems. We define the corresponding billiard map and discuss its basic properties.
We obtain a  normal form for the Hamiltonian system near the boundary $\Gamma$ and compute  the corresponding billiard map.
Dynamics near the boundary $\Gamma$ admits a version of Lazutkin's theorem \cite{Laz} on the existence of a positive measure of invariant Lagrangian tori.
As shown in \cite{Dob1,Dob2} they may serve as a basis of the short wave approximation for (\ref{eq:wave}).
We also discuss two simple examples in which the regularized flow is completely integrable.

\section{Regularization}

Let $M^n$ be a smooth Riemannian manifold with the Riemannian metric
$$
  g = \|dx\|_g^2=\langle\calA(x)dx,dx\rangle
$$
and $\phi$ a smooth function on $M$.
Then
$$
  G = \|dx\|_G^2=\frac{g}{\phi}
    = \frac{\langle\calA(x)dx,dx\rangle}{\phi(x)}
$$
 is a Riemannian metric in the domain $\Omega=\{\phi>0\}$. The metric $G$ is singular on the boundary $\Gamma=\{\phi=0\}$.
We assume that $d\phi\ne 0$ on $\Gamma$, then $\overline{\Omega}$ is a manifold with smooth boundary $\Gamma$.
The operator
$$
  \nabla_g\cdot (\phi\nabla_g u)=\Delta_Gu+\frac{n}2\langle \nabla_g\phi,\nabla_g u\rangle,
$$
in (\ref{eq:wave}) is the Laplace-Beltrami operator of the metric $G$ up to lower order terms.

The geodesic flow of the metric $G$ in $\Omega$ is incomplete.
The geodesics  are projections to $\Omega$ of  trajectories of the Hamiltonian system with the Hamiltonian\footnote
{Usually for the Hamiltonian of the geodesic flow one takes $\sqrt{H}$ or $H/2$.
With different choices only time parametrization of solutions differs.}
\begin{equation}
\label{eq:H}
  H(x,p_x) =\phi(x) \|p_x\|_g^2
           = \phi(x)\langle\calA^{-1}(x)p_x,p_x\rangle,\qquad (x,p_x)\in T^*\Omega,
\end{equation}
and the standard symplectic structure $\omega=dp_x\wedge dx$.

Since dynamics on all energy levels $H=h>0$ is the same (up to time reparametrization),
it is standard to fix the energy level $H=1$.
Then projections   of trajectories of the Hamiltonian system
are geodesics  parameterized by the arc length/2, and the length $L(\gamma)=\int_a^b\|\dot\gamma(t)\|_G\,dt$
of a geodesic $\gamma:[a,b]\to \Omega$ is the action $\int p_x\,dx$ of the corresponding trajectory in
$T^*\Omega$.

It was shown in \cite{Naz2} that the Hamiltonian flow in $T^*\Omega$ is regularizable.
We will  use a different regularization.

\begin{prop}
There exists a smooth manifold $\hat \Omega$ and a smooth map $\pi:\hat\Omega\to M$ such that:
\begin{enumerate}
\item
$\pi(\hat\Omega)=\overline{\Omega}$, $\pi$ is a double covering over $\Omega$ and
$\pi:\hat\Gamma=\pi^{-1}(\Gamma)\to\Gamma$ is a diffeomorphism.
\item
Let $\Pi:T^*(\hat\Omega\setminus\hat\Gamma)\to T^*\Omega$ be the projection
$$
  \Pi(y,p_y)=(\pi(y),{d\pi(y)^*}^{-1}p_y).
$$
Then $\hat H=H\circ\Pi$ extends to a smooth function on $T^*\hat\Omega$.
\item
The projection $\Pi$ maps trajectories of the regularized Hamiltonian  system on $T^*\hat\Omega$
with the Hamiltonian $\hat H$ and the standard symplectic form $dp_y\wedge dy$
to trajectories of the original system in $T^*\Omega$.

\item Trajectories of the regularized Hamiltonian system cross the Poincar\'e section
\begin{equation}
\label{eq:sec}
\calS_h = \{(y,p_y)\in T^*\hat\Omega:y\in \hat \Gamma\}\cap\{\hat H=h\},\qquad h>0,
\end{equation}
transversely on the energy level $\{\hat H=h\}$.
\end{enumerate}
\end{prop}

\noindent
{\it Proof}.  The regularized configuration space is defined by
$$
\hat\Omega=\{(x,\xi)\in M\times\R: \xi^2=2\phi(x)\}.
$$
Topologically this is the double of $\overline{\Omega}$ obtained by gluing two copies of $\Omega$ along $\Gamma$.
Then $\hat\Omega$ is a smooth manifold and the projection $\pi:\hat\Omega\to \overline{\Omega}$
is a smooth map. For example, if $\Omega$ is an $n$-dimensional disk, then $\hat\Omega$ is an $n$-sphere.

To compute the regularized Hamiltonian, we identify a neighborhood $U$ of $\Gamma$ in $\overline{\Omega}$
with a neighborhood $V=\{(q,r):0\le r\le \eps(q)\}$ of $\Gamma\times \{0\}$ in $\Gamma\times [0,+\infty)$
via a diffeomorphism  $\chi:V\to U$ such that $\phi(\chi(q,r))=r$ and $\chi(q,0)=q$. Then $(q,r)\in V$
defined by $x=\chi(q,r)$ serve as coordinates in $U$. In these coordinates the metric has the form
$$
  G = \frac1r\big(A(q,r)\dot r^2+2\langle \dot q,B(q,r)\rangle \dot r+\langle C(q,r)\dot q,\dot q\rangle\big).
$$
In  the corresponding symplectic coordinates in $T^*U$ we have
$$
  \omega=dp\wedge dq+dp_r\wedge dr, \qquad (q,p)\in T^*\Gamma,
$$
and the Hamiltonian (\ref{eq:H}) takes the form
$$
  H(q,p,r,p_r) = r\big(a(q,r)p_r^2+2\langle p,b(q,r)\rangle p_r+\langle c(q,r)p,p\rangle\big)
$$
for some scalar function $a$, vector function $b$, and matrix function $c$.
We make the Levi-Civita regularization
\begin{equation}
\label{L-C}
     r = \xi^2/2,\quad
   p_r = \eta/\xi,\qquad
   dp_r\wedge dr = d\eta\wedge d\xi,
\end{equation}
and obtain the regularized Hamiltonian
$$
    \hat H(q,p,\xi,\eta)
  = \frac12 a(q,\xi^2/2)\eta^2 + \langle p,b(q,\xi^2/2)\rangle\xi \eta
                               + \frac12 \langle c(q,\xi^2/2)p,p\rangle\xi^2.
$$
Trajectories of the original system with energy $H=h>0$  which hit the boundary $\Gamma$
correspond to trajectories of the regularized Hamiltonian system crossing the Poincar\'e section
\begin{equation}
\label{eq:sec2}
  \calS_h  =  \big\{(q,p,\xi,\eta): (q,p)\in T^*\Gamma, \xi=0, \eta^2=2h/a(q,0)\big\}
      \subset \{\hat H=h\}.
\end{equation}
Its two connected components $\calS_h^\pm$ are symplectically diffeomorphic to $T^*\Gamma$.
Trajectories of the regularized system cross the section transversely: $\dot\xi|_{\calS_h}= \pm\sqrt{2ha(q,0)}$.

The regularized configuration space $\hat\Omega$ is obtained by  gluing to both boundaries of
$\big\{(q,\xi):q\in\Gamma,\; |\xi|\le  \sqrt{2\eps(q)}\big\}$ a copy of $\Omega\setminus U$ via the diffeomorphism $\chi$.
\qed

\medskip

Projections of trajectories of the regularized system to $\Omega$ look like billiard trajectories,
but all trajectories hit the boundary orthogonally with respect to the metric $g$.

\section{The billiard  map}

The billiard map is the Poincar\'e first return map for the cross-section $\calS_h$ defined in (\ref{eq:sec})
or, more explicitly, (\ref{eq:sec2}). Since the first return maps corresponding to all $h>0$ are conjugate via a homothety,
we set $h=1$ and $\calS=\calS_1$. Let $\Phi^t$ be the phase flow of the regularized system.
Since $\calS$ is a section for the regularized flow, we can define a smooth Poincar\'e map $T:U\to \calS$ on an open set
$$
  U = \{z\in \calS\mid\exists t>0: \Phi^t(z)\in \calS\}.
$$
By definition   $T(z)=\Phi^\tau(z)$, where $\tau=\tau(z)$ is the  first return time.
The manifold $\calS$ is disconnected:
$$
\calS = \calS^+\cup \calS^-,\qquad  \calS^\pm = \{(q,p,\xi,\eta) : \xi=0,\; \eta =\pm \sqrt{2/a(q,0)}\}\cong T^*\Gamma.
$$
The Poincar\'e map takes  $\calS^\pm$ to $\calS^\mp$, and
 both maps $T|_{\calS^\pm}$ define the same symplectic map
$T:U\subset T^*\Gamma\to T^*\Gamma$, exactly as for the usual Birkhoff billiard.
Using the analogy with the billiard, we call $T$ the billiard map.
It is symplectic with respect to the standard symplectic form on $T^*\Gamma$.

In contrast to the usual billiard, no convexity assumptions on the boundary are needed:
the Poincar\'e map is always smooth. In fact $\Omega$ is convex with respect to the metric $G$.

Suppose the closure $\overline{\Omega}$ is compact, then $\hat\Omega$ is compact, and  the regularized flow is complete.
Then for small $\eps>0$ the set $\cup_{-\eps< t<\eps} \Phi^t(\calS)$ is open.
The Poincar\'e recurrence theorem implies that the map $T$ is defined almost everywhere on $\calS$,
so $U$ is an open dense set.

If for some $z\in\calS$ the $\omega$-limit set of $\Phi^t(z)$ is contained in $T^*\Omega$,
then $T(z)$ is undefined at $z$. If there are no compact invariant sets for the geodesic
flow of the metric $G$ in $\Omega$, then $U=\calS$ and $T:T^*\Gamma\to T^*\Gamma$ is everywhere defined.

A sufficient condition for $T$ to be globally defined is that $\Omega$ is homeomorphic to a disk and $G$ has negative sectional curvature:
in this case any geodesic of $G$ starts and ends on $\Gamma$.

\medskip

\noindent
{\bf Example}. Suppose $n=2$ and $g=|dx|^2$ is Euclidean.
Then (see e.g.\ \cite{Carmo}) the Gauss curvature of $G=g/\phi$  is given by
$$
  K =\frac{\phi}2 \Delta\ln\phi=\frac12 \Big( \Delta\phi - \frac{|\nabla\phi|^2}{\phi} \Big).
$$

Suppose   $\phi$ is concave (i.e., $-\phi$ is convex). Then $\Omega$ is a disk and $\Delta\phi < 0$. This implies $K<0$.
Hence the billiard map is everywhere defined on the cylinder $T^*\Gamma\cong S^1\times\R$, and
  $T$ is a twist map. This means that the map $T:(q,p)\mapsto (q^+,p^+)$  satisfies
$\partial q^+ / \partial p > 0$.

\smallskip

Twist symplectic self-maps of a cylinder form a popular object of studying.
In particular, they always have an infinite number of periodic \cite{Birk} and quasi-periodic \cite{Mather} trajectories.
\smallskip

Also for arbitrary $n$, if  the Hessian $h_\phi$ is negative definite,
the billiard map $T:T^*\Gamma\to T^*\Gamma$ is everywhere defined.
Indeed, the sectional curvature of $G$ corresponding to orthonormal vectors $v,w$ is given by \cite{Carmo}
$$
K_{u,v}=\frac12h_\phi(v)+\frac12h_\phi(w)-\frac{\langle\nabla\phi,v\rangle^2}{2\phi^2}
-\frac{\langle\nabla\phi,w\rangle^2}{2\phi^2}-\frac{|\nabla\phi|^2}{4\phi}<0.
$$

In section \ref{sec:average} we present  approximate formulas for the billiard map.

\section{Semigeodesic coordinates}

In this section we introduce convenient coordinates near the boundary $\Gamma$
which are similar to semi-geodesic coordinates \cite{Carmo}.

Let $L(\gamma)$ be the length of a curve $\gamma$ in $\Omega$ defined by the metric $G$
and $\rho$  the corresponding distance in $\Omega$. Let
$$
  \rho(x) = \rho(x,\Gamma)
          = \inf \{L(\gamma):\gamma\; \mbox{joins $x\in\Omega$ with $\Gamma$ in $\overline{\Omega}$}\}
$$
be the distance to $\Gamma$. It is easy to see that $\rho(x)$ is finite and
$\rho(x)\to 0$ as $x\to\Gamma$.

\begin{prop}
There exists a neighborhood $W$ of $\Gamma$ in $\overline{\Omega}$ such that:
\begin{enumerate}
\item
For any $x\in W\setminus\Gamma$ there exists a unique shortest geodesic $\gamma_x$ of the metric $G$
of length $L(\gamma_x)=\rho(x)$ joining $x$ with $q(x)\in \Gamma$.

\item
The function $s=\rho^2/4$  is smooth in $W$  and $ds\ne 0$ in $W$.

\item
The map $x\in W\setminus \Gamma\to (q(x),s(x))\in  \Gamma\times(0,+\infty)$ is a smooth diffeomorphism
on its image $V=\{(q,s):q\in\Gamma, 0<s\le \eps(q)\}$
which extends to a  smooth diffeomorphism $W\to \overline{V}=V\cup (\Gamma\times\{0\})$.
Its inverse is a  diffeomorphism $f : \overline{V} \to W$.

\item
In the coordinates $(q,s)\in\overline{V}$ on $W$ given by $x=f(q,s)$, we have
\begin{equation}
\label{eq:g/s}
  G=\frac{ds^2+\langle B(q,s)dq,dq\rangle}{s}=\frac{\tilde g}{s},
\end{equation}
where $\tilde g=ds^2+\langle B(q,s)dq,dq\rangle$ is a smooth Riemannian metric in $\overline{V}$.
\end{enumerate}
\end{prop}

\noindent{\it Proof.} We use the coordinates introduced in (\ref{L-C}). Let
$$
  \Sigma_0 = \big\{(q,p,\xi,\eta) : q\in\Gamma,\;p=0,\;\xi=0,\;\eta=\eta_0(q)\big\},\qquad \eta_0(q)=\sqrt{2/a(q,0)}.
$$
The manifold $\Sigma_0$ is isotropic: $\omega|_{\Sigma_0}=0$ and $\Sigma_0\subset\{\hat H=1\}$.
Let $\Phi^t$ be the flow of the regularized Hamiltonian system. Take small $\tau>0$ and let
$$
\Sigma=\cup_{|t|\le \tau}\Phi^t\Sigma_0.
$$
Since $\dot\xi|_{\Sigma_0}=\sqrt{2a(0,q)}>0$, the Hamiltonian vector field is transverse to $\Sigma_0$.
Hence  $\Sigma$ is a smooth Lagrangian manifold and the projection
$$
  (q,p,\xi,\eta)\in\Sigma\mapsto (q,\xi)\in\Gamma\times\R
$$
is a diffeomorphism in a neighborhood of $\Sigma_0$. Hence  a neighborhood of $\Sigma_0$ in $\Sigma$
is a smooth Lagrangian graph:
$$
          \Sigma
  \supset \big\{(q,p,\xi,\eta): p=S_q(q,\xi),\;
                   \eta=S_\xi(q,\xi),\; q\in\Gamma,\; |\xi|\le \sigma(q)\big\},
$$
where  $\sigma(q)>0$ is a smooth positive function.
The smooth generating function $S(q,\xi)$ is the action $\int p\,dq+\eta\,d\xi$ of the trajectory $\Phi^t(z_0)$
of the Hamiltonian system in $\Sigma$
starting at the point
$$
  z_0=(q_0,0,0,\eta_0)\in \Sigma_0,\qquad q_0=q_0(q,\xi),\quad  \eta_0=\eta_0(q),
$$
and ending at the point $(q,p,\xi,\eta)\in\Sigma$.  Hence for $0<|\xi|\le \sigma(q)$,
$S(q,\xi)=\rho(x)$ is the length of the  geodesic $\gamma_x$ joining the point $q_0(q,\xi)\in\Gamma$
with the point $x=f(q,\xi^2/2)\in U$.

For $\xi<0$, we have $S(q,\xi)=-\rho(x)$.
Indeed, the map  $(q,p,\xi,\eta)\mapsto (q,-p,-\xi,\eta)$ is an antisymplectic transformation conserving the Hamiltonian $\hat H$.
Hence $S(q,-\xi)=-S(q,\xi)$. It follows that $S^2(q,\xi)$ is a smooth function of $q,r=\xi^2/2$,
and so $s(x)=S^2(f^{-1}(x))/4$ is a smooth function.
Since $\partial_\xi S(q,0)=\eta_0(q)>0$, the map $(q,r)\mapsto (q_0(q,\sqrt{2r}),S^2(q,\sqrt{2r}))$
is a diffeomorphism of a neighborhood $\overline{V}$ of $\Gamma\times\{0\}$ in $\Gamma\times[0,+\infty)$.
The first two items of the proposition are proved with $\overline{W}=\chi^{-1}(\overline{V})$.

Since  $\nabla S(x)$ points inside $\Omega$ for $x\in\Gamma$, the function $\phi/s$ is smooth and positive on $W$,
and so $\tilde g=sG$ extends to a smooth Riemannian metric in $W$.
The geodesics $\gamma_x$ cross the hypersurfaces $\Gamma_s=f(\Gamma\times\{s\})$  orthogonally.
Hence in the variables $(q,s)$ introduced by the map $x=f(q,s)$, the metric $\tilde g$ is given by
$$
  \tilde g = s\,d\rho^2+\langle B(q,s)dq,dq\rangle= ds^2+\langle B(q,s)dq,dq\rangle,
$$
where $B(q,s):T_q\Gamma\to T_q^*\Gamma$ is a positive definite quadratic form smoothly depending on $(q,s)$.
\qed

\medskip

In the future we assume for simplicity that $\Gamma$ is compact. Then we can assume that $f$ is a diffeomorphism
of $\Gamma\times [0,\eps]$ onto $\overline{W}$.
Let us give a coordinate-free definition of the Riemannian metric  $g_0=\langle B(q,0)dq,dq\rangle$ on $\Gamma$.

Let $\|\cdot\|=\|\cdot\|_g$ be the norm in $TM$ defined by the metric $g$ and $\nabla=\nabla_g$ the corresponding gradient.

\begin{prop}
$$
g_0=\frac{g|_\Gamma}{\|\nabla\phi\|_g^2}=\frac{\|dq\|_g^2}{\|\nabla \phi(q)\|_g^2}.
$$
\end{prop}

\noindent
{\it Proof.}
We compute the metric $\tilde g=\frac{s}{\phi}g$ in the coordinates  defined by the map $x=f(q,s)$:
\begin{eqnarray*}
     f^*\tilde g
 &=& s\frac{\|df(q,s)\|^2}{\phi(f(q,s))}
  =  s\frac{\|(\partial_q f(q,0)+ O(s))\,dq+(\partial_s f(q,0)+O(s))\,ds \|^2}
           {\langle \nabla\phi(f(q,0)), s \partial_s f(q,0)\rangle+O(s^2)}  \\
 &=& \frac{\|\partial_q f(q,0)\,dq\|^2+\|\partial_s f(q,0)\|^2\,ds^2}
          {\langle \nabla\phi(f(q,0)),\partial_s f(q,0)\rangle} + O(s).
\end{eqnarray*}
Here $f(q,0)=q\in\Gamma$. Since $\partial_s f(q,0)\perp_g  T_q\Gamma$,
 we have  $\partial_sf(q,0)=\lambda\nabla \phi(q)$. Then
$$
    f^*\tilde g
  = \frac{\|dq\|_g^2 + \lambda^2\|\nabla \phi(q)\|^2\,ds^2}
         {\lambda\|\nabla\phi(q)\|^2} + O(s).
$$
Equation  (\ref{eq:g/s}) implies $\lambda=1$. \qed

\medskip

Geodesics of the metric $G$ in $W$, parameterized by the arc length, satisfy $\frac{d^2}{dt^2}\ln s=-s^{-1}+O(1)$.
Hence if $\delta>0$ is sufficiently small, for  $\sigma\in (0,\delta)$  the boundary $\Gamma_\sigma$
of the set $D_\sigma=\{s\ge\sigma\}$ is convex with respect to the metric $G$.

\begin{prop}
\label{prop:per}
If $\overline{\Omega}$ is compact, there exists a $2$-periodic billiard-like trajectory.
\end{prop}

\noindent
{\it Proof}.
For   $0<\sigma<\delta$, the boundary $\Gamma_\sigma$ of the set $D_\sigma$
is convex with respect to the metric $G$.  By standard arguments one can prove the existence
of a geodesic $\gamma:[0,1]\to D$ which starts and ends orthogonally (in  the metric $G$) to $\Gamma_\delta$.
Indeed, the space $X\subset W^{1,2}([0,1],D_\sigma)$ of curves  $\gamma:[0,1]\to D_\sigma$ with endpoints on $\Gamma_\delta$
is not contractible to the space of curves $\gamma:[0,1]\to \Gamma_\delta$.
By convexity, the negative gradient semilflow of the action functional $A(\gamma)=\int_0^1 \|\dot\gamma(t)\|_G^2\,dt$
is well defined and the space $X$ is positively invariant under the semiflow.
Hence there is a nontrivial critical point $\gamma$ of $A$ which is a geodesic of the metric $G$
in $D_\delta$ orthogonal to $\Gamma_\delta$ at the end points  $x_\pm\in\Gamma_\delta$.
It can be extended to the boundary via the geodesics $\gamma_{x_\pm}$ joining $x_\pm$ with $q_\pm\in\Gamma$,
and gives a periodic orbit of the regularized system,
i.e.\ a 2-periodic orbit of the billiard map: $T(q_\pm,0)=(q_\mp,0)$.
\qed

\section{Trajectories near the boundary}
\label{sec:average}

\subsection{Regularization revisited}

Let $p\in T_q^*\Gamma$, $p_s\in\R$ be the momenta conjugate to $q\in\Gamma$, $s\in[0,\eps]$.
Then the Hamiltonian corresponding to the metric $G$ in (\ref{eq:g/s}) is
\begin{equation}
\label{eq:Hs}
  H(q,p,s,p_s) = s \big( p_s^2 + K(q,p,s)\big),\qquad
  \omega = dp\wedge dq+dp_s\wedge ds,
\end{equation}
where
$$
  K(q,p,s)=\langle C(q,s)p,p\rangle,\quad C=B^{-1}.
$$

Again we perform the Levi-Civita regularization   $s=\xi^2/2$,   $p_s=\eta/\xi$, where $|\xi|\le\sqrt{2\eps}$.
Note that $|\xi|=\rho\sqrt2$, where $\rho$ is the distance  to $\Gamma$ in the metric $G$.
In the new variables we have\footnote{For simplicity we write $\hat H=H$.}
\begin{equation}
\label{eq:Hreg}
H(q,p,\xi,\eta)=\frac{\eta^2}{2}+\frac{\xi^2}{2}K(q,p,\xi^2/2),\qquad \omega=dp\wedge dq +d\eta\wedge d\xi.
\end{equation}

\subsection{Example}

As the first approximation, consider the system with Hamiltonian (\ref{eq:Hs}), where $C$ is a constant matrix.
The Hamiltonian equations have the form
$$
  \dot q = 2sCp, \quad
  \dot p = 0, \quad
  \dot s = 2sp_s, \quad
  \dot p_s = -p_s^2 - \langle Cp,p\rangle.
$$
Hence $p=\const$.  If $p\ne 0$ and $p_s(0)=0$, we obtain
\begin{equation}
\label{cycloid}
  q = q(0) + \frac{h}{k^2} \Big( t + \frac{\sin 2kt}{2k} \Big) C p, \quad
s = \frac{h}{k^2} \cos^2 kt, \quad
  p_s = -k \tan kt,\qquad |t|<\frac{\pi}{2k},
\end{equation}
where $k = \sqrt{\langle Cp,p\rangle}$  and   $H=h>0$ is the energy.
Projections of these curves to the $(q,s)$-space are cycloids.
The corresponding billiard  map $(q(-\frac{\pi}{2k}),p(-\frac{\pi}{2k}))\mapsto (q(\frac{\pi}{2k}),p(\frac{\pi}{2k}))$ is
\begin{equation}
\label{eq:T}
T(q,p)=(q+\pi  h \langle Cp,p\rangle^{-3/2}Cp,p)\qquad p\ne 0.
\end{equation}

\subsection{Action-angle coordinates}

Equations (\ref{cycloid}) determine, up to small errors, geodesics of the   metric $G$ which stay
in a small neighborhood of $\Gamma$ for a long time. Next we study these geodesics with higher precision.
Since we fixed the energy $H=h>0$, trajectories staying near the boundary have large $\|p\|$.

Let $C_0(q)=C(q,0)$, $C_1(q)=\partial_s C(q,0)$, and
$$
  K_0(q,p)=\langle C_0(q)p,p\rangle,\quad K_1(q,p)=\langle C_1(q)p,p\rangle.
$$

Up to now our equations were coordinate independent. Now we consider $q\in\R^{n-1}$ as a local coordinate in $\Gamma$,
then $p\in\R^{n-1}$.
We would like to introduce the action-angle variables $J,\ph$ in the $\xi,\eta$ plane:
$$
  \xi = \sqrt{\frac{2J}{K_0^{1/2}}}\sin\ph,\quad
  \eta = \sqrt{2JK_0^{1/2}}\cos\ph .
$$
To extend this change to a symplectic transformation  of the whole phase space,
we use the generating function
$$
    S(q,P,\ph,\eta)
  = \frac{\eta^2 \tan\ph}{2} K_0^{-1/2}(q,P) - \langle q,P\rangle.
$$
The generating function
defines the symplectic transformation $(q,p,\xi,\eta)\mapsto (Q,P,\ph,J)$ via
$$
  - p\,dq + \xi\,d\eta + J\,d\ph - Q\,dP = dS.
$$
Then
\begin{eqnarray*}
&\displaystyle
   J = S_\ph = \frac{\eta^2}{2\cos^2\ph} K_0^{-1/2}(q,P),\quad
   p = - S_q = P - \frac{\eta^2\tan\ph}{2}\frac{\partial_q K_0(q,P)}{2 K_0^{3/2}(q,P)}, & \\
&\displaystyle
 \xi = S_\eta = \eta\tan\ph K_0^{-1/2}(q,P),\quad
   Q = - S_P = q - \frac{\eta^2\tan\ph}{2}\frac{\partial_P K_0(q,P)}{2 K_0^{3/2}(q,P)}. &
\end{eqnarray*}
Although $S$ is singular at $\ph=\pm\pi/2$, the change is smooth except at $P=0$:
\begin{eqnarray*}
&\displaystyle
  \xi = \frac{\sqrt{2J}\sin\ph}{K_0^{1/4}(q,P)}, \quad
 \eta = \sqrt{2J}K_0^{1/4}(q,P)\cos\ph,  \\
&\displaystyle
    p = P - \frac{J\sin 2\ph}{4}\frac{\partial_q K_0(q,P)}{K_0(q,P)}, \quad
    Q = q - \frac{J\sin 2\ph}{4}\frac{\partial_P K_0(q,P)}{K_0(q,P)}.
\end{eqnarray*}
On the cross section $\xi=0$ we have: $P=p$, $Q=q$.

The transformation is invertible provided $J/\|P\|$ is small enough.
Note that on the fixed energy level $H=h>0$, we have $J\sim h\|P\|^{-1}$.
To obtain explicit formulas for the change, we solve the last equation for $q$ assuming that $J$ is small and $\|P\|$ large:
$$
q
  = Q + \frac{J\sin 2\ph}{4}\frac{\partial_P K_0(Q,P)}{K_0(Q,P)}
      + O\Big(\frac{J^2}{\|P\|^2}\Big),\quad
   p
  = P - \frac{J\sin 2\ph}{4}\frac{\partial_Q K_0(Q,P)}{K_0(Q,P)}
      + O\Big(\frac{J^2}{\|P\|}\Big).
$$

\subsection{Regularized Hamiltonian}

Let us rewrite the regularized Hamiltonian (\ref{eq:Hreg}) in the new variables:
\begin{eqnarray*}
      H
  &=& \frac{\eta^2}{2}+\frac{\xi^2}{2}K_0(q,P)
     + \frac{\xi^2}{2} \left( K(q,p,\xi^2/2) - K(q,P,0) \right)\\
  &=& J K_0^{1/2}(q,P) + \frac{J\sin^2\ph}{K_0^{1/2}(q,P)}
        \left( K(q,p,\xi^2/2) - K(q,P,0) \right).
\end{eqnarray*}
Let us find the main term and the first approximation:
\begin{eqnarray*}
     K_0^{1/2}(q,P)
 &=& \left(K_0(Q,P)
   + \Big\langle \partial_{Q}K_0(Q,P),\frac{J\sin 2\ph}{4}
       \frac{\partial_{P}K_0(Q,P)}{K_0(Q,P)} \Big\rangle
   + O(J^2)
     \right)^{1/2}\\
 &=& K_0^{1/2}(Q,P)
    + \frac{J\sin 2\ph}{8}
      \frac{\big\langle\partial_{Q}K_0(Q,P),\partial_{P}K_0(Q,P)\big\rangle}{K_0^{3/2}(Q,P)}
    + O\Big(\frac{J^2}{\|P\|}\Big).
\end{eqnarray*}
Next
\begin{eqnarray*}
  K(q,p,\xi^2/2) - K(q,P,0)
 &=& J\sin^2\ph \frac{\langle C_1(Q) P,P\rangle}{K_0^{1/2}(Q,P)}
    - \frac{J\sin 2\ph}{2}
      \Big\langle C_0(Q)\frac{\partial_Q K_0(Q,P)}{K_0(Q,P)},P\Big\rangle
    + O(J^2)\\
 &=& J\sin^2\ph \frac{K_1(Q,P)}{K_0^{1/2}(Q,P)}
   - \frac{J\sin 2\ph}{4}
       \frac{\langle\partial_Q K_0(Q,P),\partial_P K_0(Q,P)\rangle}{K_0(Q,P)}
   + O(J^2).
\end{eqnarray*}
Finally we obtain
\begin{eqnarray}
\nonumber
\!\! H
 &=& J K_0^{1/2}(Q,P) + J^2\sin^4\ph\frac{K_1(Q,P)}{K_0(Q,P)} \\
\label{first}
 && + \frac{J^2\sin 4\ph}{16}
        \frac{\langle\partial_Q K_0(Q,P),\partial_P K_0(Q,P) \rangle}{K_0^{3/2}(Q,P)}
    + O\Big(\frac{J^3}{\|P\|}\Big).
\end{eqnarray}

\subsection{Averaging of the regularized Hamiltonian}

On a fixed energy level $H=h>0$, we have $J\sim h\|P\|^{-1}$.
Hence $\|P\|^{-2}$ plays the role of a small parameter in the one-frequency slow-fast system with Hamiltonian (\ref{first}).
Here the angular variable $\ph$ is fast, while $Q$, $P/\|P\|$, and $J \|P\|$ are slow variables. According to the averaging principle,
one needs to replace the Hamiltonian by  $\frac1{\pi} \int_0^{\pi} H\, d\ph$ and drop higher order terms in $\|P\|^{-1}$. More precisely, we have:

\begin{prop}
\label{prop:aver}
For small $\eps>0$ and  $\|P\|>\eps^{-1}$ the system with Hamiltonian (\ref{first}) is approximated by the system with the Hamiltonian
$$
    \tilde H
  = J K_0^{1/2}(Q,P) + \frac{3J^2}8 \frac{K_1(Q,P)}{K_0(Q,P)} .
$$
This means the following. Let $Q(t),P(t),\ph(t),J(t)$ and $\tilde Q(t),\tilde P(t),\tilde\ph(t),\tilde J(t)$ be solutions of systems with Hamiltonians $H$ and $\tilde H$ respectively with the same initial condition
$Q_0,P_0,\ph_0,J_0$. Then
$$
                  \|Q(t) - \tilde Q(t)\| < \frac{c_1}{\|P_0\|^2},  \quad
  \frac{\|P(t) - \tilde P(t)\|}{\|P_0\|} < \frac{c_1}{\|P_0\|^2},  \quad
            |J(t) - \tilde J(t)| \|P_0\| < \frac{c_1}{\|P_0\|^2}
$$
on a time interval   $|t| < c_2 \|P\|^2$.
\end{prop}

\subsection{Isoenergetic reduction}

In this section we compute approximately the billiard map.
Since time parametrization is not important, we fix the energy level $H=1$ and use $\ph$ as a new time variable.
Solving the equation $H=1$ for $J$, we obtain $-J=F(Q,P,\ph)$, where
\begin{eqnarray*}
      F &=& -K_0^{-1/2}(Q,P)+ R(Q,P,\ph)+  O(\|P\|^{-5}), \\
R &=& \sin^4\ph \frac{K_1(Q,P)}{K_0^{5/2}(Q,P)} + \frac{\sin 4\ph}{16}
        \frac{\langle\partial_Q K_0(Q,P),\partial_P K_0(Q,P)\rangle}{K_0^3(Q,P)}.
\end{eqnarray*}
Then
$$
  (P\,dQ + J\,d\ph)\big|_{\{H=1\}} = P\,dQ - F\,d\ph.
$$
Thus $F$ is the new Hamiltonian, and $\ph\bmod\pi$ is the new time. Solutions $(Q(\ph),P(\ph))$ of the system with the  Hamiltonian
$F(Q,P,\ph)$ define trajectories
$$
(Q(\ph),P(\ph),\ph,J(\ph)),\qquad  J(\ph)=-F(Q(\ph),P(\ph),\ph),
$$
of our Hamiltonian system. The time parametrization is recovered by
$$
\frac{d\ph}{dt}= K_0^{1/2}(Q,P)+O(\|P\|^{-3}).
$$

For large $\|P\|>\eps^{-1}$, the variables $Q$ and $P/\|P\|$ are slow, so the averaging principle works.
Let
$$
\overline{R}(Q,P)=\frac 1\pi \int_0^\pi R(Q,P,\ph)\,d\ph.
$$
The symplectic transformation $(Q,P)\mapsto (q,p)$ with the generating function
$$
S(Q,p,\ph)=\langle Q, p\rangle +\int_0^\ph (R(Q,p,\psi)-\overline{R}(Q,p,\psi))\,d\psi,
$$
$$
P=S_Q=p+O(\|p\|^{-3}),\quad q=S_p=Q+O(\|p\|^{-2}),
$$
transforms the Hamiltonian $F$ to the normal form
$$
\calF(q,p,\ph)=N(q,p)+O(\|p\|^{-5}),\qquad N(q,p)=-K_0^{-1/2}(q,p)+\overline{R}(q,p).
$$
For $\ph=0\bmod\pi$ the transformation is identical. We obtain:

\begin{prop}
\label{prop:aver2}
For sufficiently small $\eps>0$ the billiard map $T:U\subset T^*\Gamma\to T^*\Gamma$ is well defined on the set
$U = \{(q,p)\in T^*\Gamma : \|p\|\ge \eps^{-1}\}$
and it is  close to the $\pi$-time map $\Phi_N^\pi$ of the flow $\Phi_N^\ph$ of the system with the Hamiltonian
\begin{equation}
\label{eq:N}
    N(q,p)
  = - K_0^{-1/2}(q,p)+ \frac38\frac{K_1(q,p)}{K_0^{5/2}(q,p)}=- K_0^{-1/2}(q,p)+O(\|p\|^{-3}).
\end{equation}
More precisely,
$$
T(q,p)=\Phi_N^\pi(q,p)+\big(O(\|p\|^{-4}), O(\|p\|^{-3})\big).
$$
 \end{prop}

In the first approximation $T$ is given by (\ref{eq:T}).
 We obtain a discrete version of Proposition \ref{prop:aver}.

\begin{cor}
For $m\le c\|p\|^2$ we have
$$
T^m(q,p)=\Phi_N^{\pi m}(q,p)+\big(O(\|p\|^{-2}),O(\|p\|^{-3})\big).
$$
\end{cor}

In particular, for trajectories staying near the boundary, $N$ is an approximate adiabatic invariant:
$$
N\circ T^m=N+O(\|p\|^{-3}),\qquad m\le c\|p\|^2.
$$

Recall that that $K_0(q,p)=\|\nabla\phi(q)\|^2\|p\|^2$ is the Hamiltonian describing the geodesic flow of the Riemannian metric
$g_0=\|\nabla\phi(q)\|^{-2}\|dq\|^2$, where the norm  and the gradient are defined by the metric $g$.
The Hamiltonian vector field with Hamiltonian $N$ is
$$
v_N=\frac12 K_0^{-1/2}\left(
\Big(1-\frac{15}{8}\frac{K_1}{K_0^2}\Big)v_{K_0}+\frac{3}{4K_0}v_{K_1}
\right).
$$
For large $\|p\|\sim\eps^{-1}$ we have $K_0\approx c=\const\sim\eps^{-2}$,
and so trajectories of the billiard map   are $O(\eps^{-2})$-shadowing
trajectories of the geodesic flow in $T^*\Gamma$ with the quadratic in $p$ Hamiltonian $K_0+\frac{3}{4c}K_1$.

\section{Lagrangian tori}

\subsection{Analog of Lazutkin's theorem}

Let $n=2$, then $\Gamma$ is an oriented  closed curve. Let $L$ be the total length of $\Gamma$ in the metric $g_0$
 and $q\bmod L$  the arc length parameter on $\Gamma$:
$$
L=\oint_\Gamma\frac{\|dx\|}{\|\nabla \phi(x)\|},\qquad  q=\int \frac{\|dx\|}{\|\nabla \phi(x)\|} .
$$
The second  integral is taken along an arc  of $\Gamma$ joining a fixed point $x_0\in\Gamma$ with a given point.
Then in (\ref{eq:N}) we have
$$
  K_0=p^2,\quad  K_1=k(q)p^2,\quad N(q,p)=-|p|^{-1}+\frac{3}{8}k(q)|p|^{-3},
$$
where $k(q)$ is the curvature of $\Gamma$ with respect to the metric $\tilde g$.

For definiteness let $p>0$. The billiard map $T(q,p)=(q_+,p_+)$, $q\in\T=\R\bmod L$,
is the time-$\pi$ shift\footnote
{The role of time is played by the variable $\ph$.}
along solutions of the system with the $\pi$-periodic in $\ph$ Hamiltonian
$\calF(q,p,\ph) = N(q,p) + O(p^{-5})$. Hence $T$
 is given by
$$
q_+=q+\pi p^{-2}-\frac{9\pi}{8}k(q)p^{-4}+O(p^{-6}),\quad p_+=p- \frac{3\pi}{8}k'(q)p^{-3}+O(p^{-5}).
$$
This is an exact nearly integrable twist map, so KAM-theory applies.
 The unperturbed integrable Hamiltonian $\tilde F$ is nondegenerate.
This means that $\nu'(\eps)\ne 0$, where  $\nu(\eps)=\eps^{-2}+O(\eps^{-4})$
is the frequency on the unperturbed invariant circle $\{N=\eps\}$.

By the standard KAM-theory (see e.g.\ \cite{AKN,Moser}), we obtain:

\begin{prop}
The billiard map has a Cantor family of 1-dimensional invariant tori $\calT_\eps\subset \T\times \R$,
$$
\calT_\eps=\{(q,p):q\bmod L,\; p^{-1}=\eps+\frac{3}{8}k(q)\eps^3+ O(\eps^{5})\},
$$
where $\eps$ belongs to a Cantor set with density point $0$.
\end{prop}

The gaps in the Cantor set are exponentially small in $\eps$ as $\eps\to 0$.
This statement is an analog of the Lazutkin theorem \cite{Laz} (see also generalizations in \cite{Popov})
for traditional Birkhoff billiards  with smooth strictly convex boundary.

Any torus $\calT_\eps$ is obtained as an intersection of a 2-dimensional  invariant torus
$\calT_\eps^2\subset T^*\hat\Omega$ of  the regularized system  with the cross section $\calS$.
The  action variables $I_i=\oint_{\gamma_i}p_x\,dx$ on $\calT^2_\eps$ are
$$
  I_1=\pi\eps+O(\eps^3),\quad I_2=\frac{L}{\eps}-\frac{3\eps}{8}\int_0^L k(q)\,dq+O(\eps^3).
$$
Similar to Lazutkin's theorem, projections of trajectories in $\calT_\eps^2$ to $\Omega$, i.e.\ corresponding billiard trajectories in $\Omega$,
are tangent to the Cantor family of caustics
$$
  \{x\in\Omega:\rho(x,\Gamma)=\eps+O(\eps^2)\}.
$$

\subsection{Tori around tori}

A positive measure of Lagrangian tori appears near any non-degenerate reducible Diophantine invariant torus, see e.g.\ \cite{AKN}.
For example, if the periodic orbit from Proposition \ref{prop:per} is elliptic, then generically it is surrounded by a Cantorian family of Lagrangian tori.
Or, if there is a nondegenerate (in the KAM sense) elliptic closed geodesic of the metric $g_0$ on the boundary $\Gamma$,
then, as in \cite{Popov}, using Proposition \ref{prop:aver2},  we will obtain a Cantor family of $(n-1)$-dimensional Lagrangian tori
in $T^*\Gamma$ for the billiard map, and hence a family of $n$-dimensional invariant Lagrangian tori for the regularized flow in $T^*\hat\Omega$.

\section{Integrable cases}

There are two simple situations when the regularized system is completely integrable, both in the case $M=\R^n$.
For $n=2$ they were studied in \cite{Dob1,Dob2}.

\medskip

{\bf a}. Suppose $g$ and $\phi$ are $SO(n)$-invariant.
Then $\Omega$ is an $n$-dimensional ball with center at the origin or a domain between two concentric $(n-1)$-spheres.
The corresponding regularized geodesic flow is completely integrable due to the existence of the algebra
of linear in momentum Noether integrals.

\medskip

{\bf b}. Suppose $g=\sum g_i(x_i)\,dx_i^2$. Let $\phi_i:[-1,1]\to\R$ be smooth  functions
such that $\phi_i>0$ on $(-1,1)$, $\phi_i(\pm 1)= 0$, and $\phi'_i(\pm 1)\ne 0$. Set $\chi_i=1/\phi_i$ and
$$
\phi(x_1,\dots,x_n) = (\chi_1(x_1)+ \dots + \chi_n(x_n))^{-1}.
$$
Then $\phi>0$ on $\Omega=(-1,1)^n$ and $\phi|_\Gamma=0$. We have
$$
  G =\frac{g}{\phi}=(\chi_1(x_1)+ \dots + \chi_n(x_n)) (g_{1}(x_1)\, dx_1^2+\dots+g_n(x_n)\,dx_n^2) .
$$
In the corresponding Hamiltonian the variables separate. Hence the regularized geodesic flow is completely integrable.
In this case $\Omega = (-1,1)^n$ is an $n$-dimensional cube.
The billiard map is not  defined at the corner singularities, although  they are probably regularizable.

\medskip

Other complete integrability cases are obtained taking direct products of these systems.
It will be interesting to find less evident integrable examples.

\section*{Funding}

The work of D.Treschev was supported by the Russian Science Foundation under grant no. 20-11-20141,\newline
https://rscf.ru/en/project/20-11-20141/. The work of S.Bolotin  was performed at the Steklov International Mathematical Center and
supported by the Ministry of Science and Higher Education of the Russian Federation (agreement no. 07515-2022-265).
Sections 1, 2, 4, 5.6, 6.1, 6.2 were written by S.~Bolotin, and sections  3, 5.1, 5.2, 5.3, 5.4, 5.5, 7 by D.~Treschev.

\end{document}